\documentclass[12pt]{article}
\usepackage{amsmath}
\usepackage{amscd}
\usepackage{amssymb}
\usepackage{array}

\newtheorem{theorem}{Theorem}[section]
\newtheorem{example}{Example}[section]
\newtheorem{proposition}{Proposition}[section]
\newtheorem{lemma}{Lemma}[section]
\newtheorem{remark}{Remark}[section]

\begin{document}

\newcommand{\R}{\mathbb{R}}
\newcommand{\C}{\mathbb{C}}
\newcommand{\Z}{\mathbb{Z}}

\newcommand{\g}{\mathcal{G}}
\newcommand{\A}{\mathcal{A}}
\newcommand{\B}{\mathcal{B}}
\newcommand{\F}{\mathcal{F}}
\newcommand{\I}{\mathcal{I}}
\newcommand{\Ha}{\mathcal{H}}
\newcommand{\Oc}{\mathcal{O}}
\newcommand{\m}{\mathbf{m}}
\newcommand{\La}{\mathcal{L}}
\newcommand{\N}{\mathcal{N}_{\Theta^r}}
\newcommand{\U}{\mathcal{U}}

\newcommand{\sym}{\mbox{sym}}
\vskip 2cm
  \centerline{\LARGE \bf Deformation Quantization of non Regular }\bigskip
  \centerline
{ \LARGE\bf  Orbits of Compact Lie Groups}

\vskip 2.5cm

\centerline{ M. A. Lled\'o}

\bigskip

\centerline{\it Dipartimento di Fisica, Politecnico di Torino,}
\centerline{\it Corso Duca degli Abruzzi 24, I-10129 Torino,
Italy, and} \centerline{\it INFN, Sezione di Torino, Italy.}
\centerline{{\footnotesize e-mail: lledo@athena.polito.it}}
 \vskip 3cm

\begin{abstract}

In this paper we construct a deformation quantization of the
algebra of polynomials of an arbitrary (regular and non regular)
coadjoint orbit of a compact semisimple Lie group.  The deformed
algebra is given as a quotient of the enveloping algebra by a
suitable ideal.
\end{abstract}

\vfill\eject

\section{Introduction}

Coadjoint orbits of Lie groups can model phase spaces of physical
systems since they are symplectic manifolds with a Hamiltonian
symmetry group, the original Lie group $G$ itself. Not only do
they enjoy nice properties at the classical level, but they are
extremely interesting as quantum systems.   The Kirillov-Kostant
principle associates in many cases a unitary representation of $G$
to the orbit, and this representation is the starting point of the
quantization of the system. The algebra of quantum observables is
related to the enveloping algebra $U$ of of the Lie algebra $\g$
of $G$. In fact, it can be expressed as a quotient of $U$ by a
suitable ideal $\I$ \cite{vo}. This ideal belongs to the kernel of
the associated representation.

Given such  algebra of observables, one may wonder if the approach
of  deformation quantization, which does not make reference to the
Hilbert space, can give an isomorphic algebra. There are at least
two paths  that one can follow at this point. One is to make a
differential deformation in the sense of Bayen, Flato, Fronsdal,
Lichnerowicz and Sternheimer \cite{bffls}, De Wilde and Lecomte
\cite{dl}, Fedosov \cite{fe} and Kontsevich \cite{ko}. Progress
were made in the study of tangential differential deformations in
\cite{fr, cgr, alm, cgr2, ma}. The other approach is to see the
algebra $U/\I$ as some kind of deformation of the polynomials on
an algebraic manifold (the coadjoint orbit itself for the case of
compact groups).

In \cite{fl} the second approach was taken and a family
deformations  was constructed in this way for regular orbits of
compact Lie groups. The family contained as a particular case the
star product found  in Ref.\cite{cg}. It was also shown in
\cite{cg} that such  star product is not differential, so it
seemed that there was some kind of incompatibility among the two
approaches, the differential  and the algebraic one, the last one
making use of the structure of the enveloping algebra. These star
products were further studied in Ref.\cite{fll}, and the problem
of the compatibility or equivalence among these two kinds of
deformations is studied in Ref.\cite{fl2}. There, it was proven
that there is an injective homomorphism from one algebraic  star
product into the differential star product algebra.

All these constructions rely heavily on the regularity of the
orbit. Nevertheless, the structure of the coadjoint representation
is much richer.  There are many  interesting cases of symplectic
manifolds with symmetries which are  diffeomorphic to non regular
orbits (or non generic orbits, with dimension less than the
maximal one). So it is of great interest to see their quantization
in some way. Given a star product, one can define the star
exponential and then a star representation of the original group.
Much work has been done in star representation theory of
semisimple groups (See for example Refs.\cite{acg, mo,mo2, cn}).

Also, some precursor of this result may be found in Ref.\cite{hk}
and for the differential version in Refs.\cite{co,na}.

In this paper we solve this problem by generalizing the algebraic
approach of \cite{fl,fll} to the case of non regular orbits. In
Section 2 we make a short summary of the properties of the
coadjoint orbits and present the non regular orbits in the
appropriate way for our purposes. As an example, we give the
coadjoint orbits of SU$(n)$. In Section 3 the existence of the
deformation quantization is proven by showing that a certain
algebra $U_h/\I_h$ has the right properties, in particular the one
of being isomorphic as a $\C[h]$-module to the polynomials on the
orbit.

Finally we want to note that the construction could be extended to
semisimple orbits of non compact groups, by quantizing the real
form of the complex orbit, which is in fact a union of connected
components each of them a real orbit. The case of nilpotent orbits
is for the moment not considered in this approach.

\section{Regular and non regular orbits as algebraic varieties}
Let $G$ be a compact semisimple group of dimension $n$ and rank
$m$ and $\g$ its Lie algebra. Let $\g^*$ be the dual space to
$\g$. On $C^\infty(\g^*)$ we have the Kirillov Poisson structure
$$
\{f_1,f_2\}(\lambda)=\langle[(df_1)_{\lambda},(df_2)_{\lambda}],\lambda\rangle,
\qquad f_1,f_2 \in C^{\infty}(\g^*), \quad \lambda \in \g^*. $$
$df_{\lambda}:\g^*\rightarrow \R$ can be considered as an element
of $\g$, and $[\; ,\,]$ is the Lie bracket on $\g$. Let
$\{X_1\dots
 X_n\}$ be   a basis of $\g$ and $\{x^1,\dots x^m\}$ the
 coordinates on $\g^*$ in the dual basis. We have that
 $$\{f_1,f_2\}(x^1,\dots x^n)=\sum_{ijk}c_{ij}^k x^k\frac{\partial f_1}
 {\partial x^i }\frac{\partial f_2}{\partial x^j
 },\qquad{\rm where}\qquad
 [X_i,X_j]=\sum_kc_{ij}^kX_k.$$
 Notice that the ring of polynomials of $\g^*$, ${\rm Pol}(\g^*)$, is closed under
 the  Poisson bracket.

 The Kirillov Poisson structure is not symplectic nor regular. As
 any Poisson manifold, $\g^*$ can be foliated in symplectic
 leaves, the Poisson bracket restricting to a symplectic Poisson
 bracket on the leaves of the foliation. The symplectic leaves
 support a hamiltonian action of $G$. Indeed, they are the orbits of
 the coadjoint action of $G$ on $\g^*$.

The coadjoint action of $G$ on $\g^*$, ${\rm Ad}^*$, is defined by
$$ \langle{\rm Ad}^*_g\lambda ,Y\rangle=\langle\lambda,{\rm
Ad}_{g^{-1}}Y\rangle \quad \forall\; g\in G,\quad \lambda\in
\g^*,\quad Y\in \g. $$ We denote by $C_\lambda$ the orbit of an
element $\lambda\in\g^*$ under the coadjoint action. The coadjoint
orbits are real irreducible algebraic varieties (see for example
\cite{bo, va1}). Let ${\rm{Inv}}(\g^*)\subset {\rm{Pol}}(\g^*)$
the subalgebra of
 polynomials on $\g^*$ invariant under the coadjoint action.
 Then ${\rm{Inv}}(\g^*)=\R[p_1,\dots p_m]$, where $\{p_1\dots p_m\}$
 is a system of algebraically independent homogeneous invariant
 polynomials, and $m$ is  the rank of $\g$ (Chevalley's theorem).

  Since
 $G$ is a semisimple Lie group, we can identify $\g\simeq \g^*$ by
 means of the invariant Cartan-Killing form, the isomorphism intertwining
 the
 adjoint and coadjoint representations. From now on we will assume that such identification has been made.
The set ${\rm{Inv}}(\g^*)$ is in one to one correspondence with
the set of polynomials on the Cartan subalgebra that are invariant
under the Weyl group, the isomorphism being given by the
restriction.

Consider the adjoint action of $\g$ on $\g$, ${\rm ad}$. The
characteristic polynomial of $X\in \g$ in the indeterminate $t$ is
 $$\det(t\cdot {\bf 1}-{\rm{ad}}_X))=\sum_{i \geq m}
q_i(X)t^i. $$  The $q_i$'s are invariant polynomials. An element
$X\in \g$ is regular if $q_m(X) \neq 0$. When restricting to the
Cartan subalgebra $\Ha$, $q_m$ is of the form

$$q_m(H)=\prod_{\alpha\in \Delta}\alpha(H), \qquad H\in \Ha,$$
with $\Delta $ the set of roots. It is clear that an element $H$
is regular if and only if it belongs to the interior of a Weyl
chamber. Any orbit
 of the adjoint action intersects the Weyl chamber in one and only one point.
 If $H$ is regular, ${\rm Ad}_g=gHg^{-1}$ is also regular, so the
 orbit $C_H$ is a regular orbit. The differentials of Chevalley's generators $dp_i$, $i=1,\dots m$ are
linearly independent only on  the regular elements \cite{va2}.
Moreover, the regular orbits are defined as algebraic varieties
 by the polynomial equations
$$p_i=c^0_i,\qquad c^0_i\in \R \qquad i=1,\dots m.$$
  The ideal
of polynomials vanishing on a regular orbit is given by
$$\I_0=(p_i-c_i,\; i=1,\dots m),$$ and the coordinate ring of
$C_X$ is ${\rm Pol}(C_X)\simeq{\rm Pol}(\g^*)/\I_0$.

\begin{example} Coadjoint orbits of SU(n)\end{example} We will
consider the compact Lie group SU(n), with (complexified) Lie
algebra $\A_m={\rm sl}(m+1, \C)$, $n=m+1$. A Cartan subalgebra of
$\A_m$ is given by the diagonal matrices $$\Ha=\{{\rm
diag}(a_1,a_2,\dots a_{m+1}), \quad a_1+a_2+\dots +a_{m+1}=0\quad
a_i\in \C\}.$$ Denoting by $$\lambda_i({\rm diag}(a_1,a_2,\dots
a_{m+1})=a_i,$$ then the roots are given by
$\alpha_{ij}=\lambda_i-\lambda_j$. We will denote the simple roots
as $\alpha_i=\lambda_i-\lambda_{i+1}$. The {\it root vectors} are
defined by $$\alpha_{ij}(H)=\langle H,H_{ij}\rangle,$$ where
$\langle\;,\;\rangle$ denotes the Cartan-Killing form on  $\Ha$.
One normalizes the root vectors as $$\bar H_{ij}=\frac{2}{\langle
H_{ij},H_{i,j}\rangle}H_{ij}.$$

The Weyl group is the group of reflections of $\Ha$ generated by
$$\omega_{ij}(H)=H-\alpha_{ij}(H)\bar H_{ij}.$$ We have
$$\omega_{ij}({\rm diag}(a_1,\dots a_i,\dots a_j ,\dots a_{m+1}))=
{\rm diag}(a_1,\dots a_j,\dots a_i ,\dots a_{m+1}),$$ so an
element of the Weyl group is $$\omega_s({\rm diag}(a_1, a_2,\dots
a_{m+1}))=({\rm diag}(a_{s^{-1}(1)}, a_{s^{-1}(2)}\dots a_i,\dots
a_{s^{-1}(m+1)})), \quad s\in \Pi_{m+1},$$ with $\Pi_{m+1}$  the
group of permutations  of order $m+1$.

We take the real span $\Ha^\R=\oplus_{i=1}^l\R\bar H_i$. Any point
$(a_1,a_2,\dots a_{m+1})\in \Ha^R$ can be brought to the form
$$(a_1,a_2,\dots a_{m+1}), \qquad (a_1\geq a_2\dots \geq
a_{m+1})$$ by applying a suitable permutation.  The intersection
of this subset of $\R^n$ with the solution of $a_1+a_2+\dots
+a_{m+1}=0$ is the closed principal Weyl chamber since
$$\alpha_i(a_1,a_2,\dots a_{m+1})\geq 0 \qquad \forall i.$$

We go now to the compact real form. By means of the map
$$A+iB\longrightarrow \begin{pmatrix}A&B\\-B&A\end{pmatrix}\in
{\rm SO(2n)}, \qquad A+iB\in {\rm SU(n)},$$ (real representation
of SU$(n)$, $n=m+1$) we obtain an isomorphism $${\rm
SU(n)}\simeq{\rm SO(2n)}\cap{\rm Sp(2n)}.$$ The Cartan subalgebra
is represented by matrices of the form $${\rm diag}
(a_1\Omega,a_2\Omega, \dots a_n\Omega), \quad \Omega
=\begin{pmatrix}0&1\\-1&0\end{pmatrix}, \quad a_1+a_2+\dots
+a_{n}=0, \quad a_i\in \R.$$ To compute the orbits of the Weyl
group it is enough to restrict to one  Weyl chamber. Consider the
partition $n=p_1+\cdots +p_r$, $p_i$ positive integers, and
consider a  point $H=(a_1,a_2,\dots a_{n})$ in the Weyl chamber
such that $$a_1=\cdots = a_{p_1}\geq a_{p_1+1}=\cdots
=a_{p_1+p_2}\geq \cdots =a_n.$$ The group of matrices leaving $H$
invariant by the adjoint action, $gHg^{-1}$ is isomorphic to ${\rm
Sp}(2p_1)\times {\rm Sp}(2p_2)\times {\rm Sp}(2p_r)\subset {\rm
Sp}(2n) $. So the isotropy group of $H$ in SU(n) is
\begin{equation}{\rm SO}(2n)\cap {\rm Sp}(2p_1)\times
{\rm Sp}(2p_r)={\rm S}({\rm U}(p_1)\times \cdots {\rm
U}(p_r)).\label{isotropy}\end{equation}

The coadjoint orbits are spaces of the type $G/H$ with $G={\rm
SU}(n)$ and $H$ one of the isotropy groups in  (\ref{isotropy}).
The regular orbits are diffeomorphic to SU(n)/U(1)$\times \cdots
\times$U(1) (with $n-1$ U(1) factors).  The non regular orbits
correspond obviously to the border of the Weyl chamber, where at
least one of the roots has value   zero, $\alpha_i(H)=0$.
\hfill$\blacksquare$

\bigskip
Non regular orbits are also algebraic varieties, but unlike the
regular orbits, the ideal $\I_0$ of a non regular orbit is not
generated by the invariant
 polynomials $(p_i-c_i,\; i=1,\dots m)$. The next proposition shows that there is a special set of
  generators that are invariant as a set.
\begin{proposition}If  $C_\lambda$ is a semisimple coadjoint  orbit of a semisimple Lie group, then
the ideal of $C_{\lambda}$, $\I_0$, is generated by polynomials $r_1(x),\dots r_l(x)$ such that
\begin{equation}r_\alpha(g\cdot x)=\sum_{\beta=1}^l T(g)_{\alpha\beta}r_\beta(
x), \qquad  \alpha=1,\dots l,\qquad g\in
G.\label{polyrep}\end{equation} where $T$ is a representation of
$G$ (and of $\g$).
\end{proposition}
{\it Proof.}  By the Hilbert basis theorem, every ideal in
$\C[x_1\dots x_n]$ has a finite generating set. Let
$\{q_1(x),\dots q_r(x)\}$ be an arbitrary finite set of generators
 of $\I_0$. Let $q(x)\in \I_0$. $G$ is an algebraic group and the action of $G$ on $\g$ sends polynomials into polynomials.
 Then $q^g(x)=q(g^{-1}x)\in \I_0$. We consider the set
$$\{q^g_i(x), \quad i=1,\dots r,\quad g\in G\}$$ which obviously
generates $\I_0$. We consider the $\C$-linear span of $\{q_i^g\}$.
Notice that it is a finite dimensional vector space since the
degree of $q_i$ doesn't change under the action of the group.
 We take  a $\C$-basis of it denoted  by $r_1(x),\dots r_l(x)$. Then
$$r_j^{g}(x)=r_j(g^{-1}x)=\sum_{k=1}^{l}T_{jk}(g^{-1})r_k(x),
\qquad j=1,\dots l. $$ It is immediate to see  that the matrices
$T(g)$ form a representation of $G$, as we wanted to
show.\hfill$\blacksquare$

\bigskip

In terms of the Lie algebra, equation (\ref{polyrep}) can be written as
\begin{equation}X.r_i(x)=\sum_kT(X)_{ik}r_k(
x),\qquad X\in \g\label{polyrep2}\end{equation} where $X$ acts as a derivation on $\C[\g^*]$.

\section{Deformation quantization of non regular orbits}

In this section we consider algebras over $\C[h]$ (polynomials on
the indeterminate $h$). We will consider the complexification of
the polynomial ring of the orbit. The deformation quantization
that we will obtain  is an algebra over $\C[h]$, so it  can be
evaluated at any particular value $h=h_0\in \R$.

 We consider the
tensor algebra $T_{\C}(\g)[h]$ and its proper two sided ideal
\begin{equation} \mathcal{L}_{h}=\sum_{X,Y \in \g}
T_{\C}(\g)[h] \otimes(X \otimes Y - Y \otimes X - h[X,Y]) \otimes
T_{\C}(\g)[h]. \label{lh}\end{equation} We define $U_{h}=
T_{\C}(\g)[h]/\mathcal{L}_{h}$. It is well known that $U_h$ is a
deformation quantization of $\C[\g^*]={\rm Pol}(\g^*)$. In
\cite{fl} it was shown that if $C_\lambda$ is a regular orbit,
then there exists a deformation quantization
 ${\rm Pol}(C_\lambda)$ as $U_h/\I_h$ where
$\I_h\rightarrow \I_0$ when $h\rightarrow 0$.  In this section we
want to
 generalize that construction to the case of non
regular orbits. First we briefly review how the deformation is
obtained in the regular case.

Let $\{X_1,\dots , X_n\}$ be a basis for $\g$ and let $\{x_1,\dots
, x_n\}$ be the corresponding generators of  $\C[\g^*]$. There is
a natural isomorphism of  $\C[\g^*]$ with the symmetric tensors
$ST_{\C}(\g)\subset T_{\C}(\g)$, $\mbox{Sym}
:\C[h][\g^*]\longrightarrow ST_{\C[h]}(\g)$  by
\begin{equation} \hbox{Sym}(x_1\cdots x_p)=\frac{1}{p!}\sum_{s \in
S_p}X_{s(1)}\otimes\cdots\otimes X_{s(p)} \label{sym}
\end{equation} where $S_p$ is the group of permutations of order
$p$. The composition of the symmetrizer with the natural
projection $T_{\C[h]}(\g)\longrightarrow U_h$ is a linear
isomorphism $W:\C[h][\g^*]\rightarrow U_h$ called the Weyl map.
This map has the following interesting property  (see for example
Ref. \cite{va3} for a proof). Let $A$ be  an automorphism
(derivation) of $\g$. It extends to an automorphism (derivation)
of $U_h$ denoted by $\tilde A$. It also extends to an automorphism
(derivation) $\bar A$ of  $ST_\C(\g)\simeq \C[\g^*]$. Then
\begin{equation}W\circ\bar A=\tilde A\circ
W.\label{intertwin}\end{equation}

Taking $A= {\rm{ad}}_X$, $X\in \g$,  (\ref{intertwin}) implies
that the images of the invariant polynomials $P_i=W(p_i)$,
$i=1,\dots, m$ belong to the center of $U_h$. The two sided ideal
$\I_h=(P_i-c_i(h))$, $c_i(h)\in \C[h]$ with $$\begin{CD}
\I_h@>>{h\mapsto 0}>\I_0\end{CD}$$
 is equal
to the right and left ideals with the same generators,
$I_h=\I_h^L=\I_h^R$.  We have that $U_h/\I_h$ is a deformation
quantization of ${\rm Pol}(C_\lambda)$ for $C_\lambda$ a regular
orbit \cite{fl, fll}. The subtle point in the proof was to show
that there is a $\C[h]$-module isomorphism $\C[\g^*][h]/\I_0\simeq
U_h/\I_h$. This was done by choosing a basis in
$\C[\g^*][h]/\I_0$, mapping it to $U_h/\I_h$ and showing that the
image is a basis of  $U_h/\I_h$. Proving the linear independence
of the images made use of the regularity hypothesis. Here we will
give a proof of the independence that does not use the regularity
condition, so it applies for non regular orbits as well.

 Let $C_\lambda$ be a non
regular orbit, and let its ideal $\I_0$ be generated by
 $r_\alpha$, $\alpha=1,\dots l$ satisfying the condition (\ref{polyrep}).  We consider
 the images under the Weyl map of the generators, $R_\alpha=W(r_\alpha)$. We have the following

\begin{lemma}
{ The left and right ideals $\I_h^{L,R}$ generated by $R_\alpha$,
$\alpha=1,\dots l$ are equal and then equal to the two sided ideal $\I_h$.}\end{lemma}

 {\it Proof.}  It is enough to prove that $[X,R_\alpha]=\sum_\beta C_{\alpha\beta}R_\beta $ for any $X\in \g$
  and some $C_{\alpha\beta}\in \C$. Using (\ref{intertwin}) and (\ref{polyrep2}) for
$A={\rm ad}_X$, we have
\begin{eqnarray*}[X, R_\alpha]=[X, W(r_\alpha)]=W(X.r_\alpha)=W(T(X)_{\alpha\beta}r_\beta)=\\
T(X)_{\alpha\beta}W(r_\beta)=
T(X)_{\alpha\beta}R_\beta.\end{eqnarray*} \hfill$\blacksquare$

\bigskip

In order to show that $U_h/\I_h$ is a deformation quantization of
${\rm{Pol}}(C_\lambda)$ we have to show that
 ${\rm Pol}(C_\lambda)[[h]]$ isomorphic as a $C[[h]]$-module to
$U_h/\I_h$. We will do it in several steps.
We need to introduce the   evaluation map ${\rm{ev}}_{h_0}:U_h\rightarrow U_h/(h-h_0)\simeq U_{h_0}$ where
 $U_{h_0}$ is the enveloping algebra of $\g$ with bracket $h_0[\;,\;]$. As it is well known for any enveloping algebra,
  $U_{h_0}$ is a filtered algebra,
$$\g\simeq U_{h_0}^{(1)}\subset U_{h_0}^{(1)}\subset \cdots
\subset U_{h_0}^{(n)}\subset \cdots$$ where
$$U_{h_0}^{(n)}/U_{h_0}^{(n-1)}\simeq ST^{(n)}_\C(\g)\simeq
\C[\g^*]^{(n)}.$$ The graded algebra associated to the filtered
algebra $U_{h_0}$ is then the algebra of symmetric tensors on $\g$
(or polynomials on $\g^*$), so there exists a natural projection
\begin{equation}\pi:U_{h_0}\rightarrow
\C[\g^*].\label{proyection}\end{equation}

\bigskip

We introduce the graded lexicographic ordering in  $\C[\g^*]$. We
consider the ideal generated by the leading terms of $r_\alpha$,
$(LT(r_\alpha))$.  Any equivalence class in $\C[\g^*]/\I_0$ has a
unique representative as a linear combination of elements  in the
set $\B=\{x_{i_1}x_{i_2}\cdots x_{i_k},\; x_{i_1}x_{i_2}\cdots
x_{i_k}\notin (LT(r_\alpha))\}$. (The basis $\{r_\alpha\}$ can be
chosen as a Groebner basis). In fact, the elements of $\B$ are
linearly independent over $\C$ and form a basis of $\C[\g^*]/\I_0$
(see for example \cite{clo} for a proof). We will denote by $J$
the set of indices $(i_1,\dots i_k)$ of elements of $\B$.

We have the following
\begin{lemma} The standard monomials
 \begin{equation}\{X_{i_1}\cdots X_{i_k}, \quad (i_1,\dots i_k)\in
 J\}\label{mono}\end{equation}
are linearly independent  in $U_{h_0}/\I_{h_0}$\label{lemmali}.
\end{lemma}

{\it Proof. } Suppose that there exists a linear relation among them
\begin{equation}\sum_{(i_1\dots i_k)\in J}c^{i_1\dots i_k}X_{i_1}\cdots X_{i_k}=\sum_\alpha A^\alpha R_\alpha,
\qquad A^\alpha\in U_{h_0}.\label{lide}\end{equation} We project
(\ref{lide}) onto $\C[\g^*]$ as in  (\ref{proyection}),
\begin{eqnarray}\pi(\sum_{(i_1\dots i_k)\in J}X_{i_1}\cdots
X_{i_k})&=& \sum_{(i_1\dots i_k)\in J_0}c^{i_1\dots
i_k}x_{i_1}\cdots x_{i_k},\nonumber\\ \pi(\sum_\alpha A^\alpha
R_\alpha)&=& \sum_{\alpha} b^{\alpha} r_\alpha + \hbox{terms with
degree } <l_0 \label{medium}
\end{eqnarray} where $J_0\subset J$ corresponds to  the monomials
with highest degree ($l_0$) in the left hand side of (\ref{lide}).
The second equation in (\ref{medium})  expresses the fact that the
right hand side of (\ref{lide}) must project to a linear
combination of $r_\alpha$ modulo terms of lower degree. So we have
\begin{equation}\sum_{(i_1\dots i_k)\in J_0}c^{i_1\dots
i_k}x_{i_1}\cdots x_{i_k}=\sum_{\alpha} b^{\alpha} r_\alpha +
\hbox{terms with degree } <l_0 \label{result}.\end{equation}
Taking  the leading term of both sides of the equation
(\ref{result}), and taking into account that $\{r_\alpha\}$ is a
Groebner basis, we obtain that the leading term of
$\sum_{(i_1\dots i_k)\in J_0}c^{i_1\dots i_k}x_{i_1}\cdots
x_{i_k}$ must be proportional to one of the leading terms
$LT(r_\alpha)$, which is not possible by the construction of the
basis $\B$.\hfill$\blacksquare$

We can prove now the independence of the monomials (\ref{mono}) on
$U_h/\I_h$.
\begin{proposition} The standard monomials
$$\{X_{i_1}\cdots X_{i_k}, \quad (i_1,\dots i_k)\in
 J\}$$
are linearly independent  in $U_{h}/\I_h$. \label{propo1}
\end{proposition}

{\it Proof}. Suppose that there is a linear combination of them
equal to zero,

\begin{equation}\sum_{(i_1\dots i_k)\in J}c^{i_1\dots i_k}(h)X_{i_1}\cdots X_{i_k}=\sum_\alpha A^\alpha R_\alpha,
\qquad A^\alpha\in U_{h}.\label{lideh}\end{equation} Applying the
evaluation map ${\rm ev}_{h_0}$ to (\ref{lideh}), we have

$$\sum_{(i_1\dots i_k)\in J}c^{i_1\dots i_k}(h_0)X_{i_1}\cdots
X_{i_k}=\sum_\alpha A_0^\alpha R_\alpha, \qquad {\rm
ev}_{h_0}A^\alpha=A_0,$$ which implies, by Lemma \ref{lemmali}
that $c^{i_1\dots i_k}(h_0)=0$. Since this is true for  infinitely
many $h_0$ and $c^{i_1\dots i_k}(h)$ is a polynomial, we have
that $c^{i_1\dots i_k}(h)=0$.

\hfill$\blacksquare$

We want now to prove that the monomials (\ref{mono}) generate
$U_h/\I_h$, so they form a basis of the $\C[h]$-module which is
then free and isomorphic to $\C[\g^*]/\I_0[h]$. The proof is the
same than the one in Ref. \cite{fl} for regular orbits, so we do
not repeat it here

\begin{proposition} The  standard monomials
$\{X_{i_1}\cdots X_{i_k}\}$ with ${(i_1,\dots , i_k) \in J}$
generate $U_h/I_h$ as
$\C[h]$-module.\label{propo2}\end{proposition}

Summing up, we have the following

\begin{theorem}
Let $\Theta$ be a (possibly non regular) coadjoint orbit  of a
compact semisimple group. In the same notation as above, $U_h/I_h$
is a deformation quantization of ${\rm
Pol(\Theta)}=\C[\g^*]/\I_0$. It has the following properties:

\smallskip

\noindent 1. $U_h/I_h$ is isomorphic to $\C[\g^*]/I_0[h]$ as a
$\C[h]$-module.

\smallskip

\noindent 2. The multiplication in $U_h/I_h$ reduces ${\rm
mod}(h)$ to the one in $\C[\g^*]/I_0$.

\smallskip

\noindent 3. The bracket $[F,G]=FG-GF$, in   $U_h/I_h$, reduces
${\rm mod} (h^2)$ to ($h$ times)  the Poisson bracket on the
orbit.

\end{theorem}

{\sl Proof.}
 1 is a consequence of Propositions \ref{propo1} and
\ref{propo2}. 2 and 3 are trivial.

\begin{remark}Extension to $\C[[h]]$\end{remark}
The extension to $\C[[h]]$ can be made by taking the inverse
limits of the families $U_h/h^nU_h$ and
$(U_h/\I_h)/h^n(U_h/\I_h)$. The elements $\{X_{i_1}\cdots
X_{i_k}\}$ with ${(i_1,\dots , i_k) \in J}$ are linearly
independent in the inverse limit since they are so in each of the
projections to $(U_h/\I_h)/h^n(U_h/\I_h)$. Then one can show that
they form a basis.

\section*{Acknowledgements}
I want to thank V. S. Varadarajan and R. Fioresi for useful
discussions.

\end{document}